\documentclass[letterpaper, 10 pt, conference]{ieeeconf} 
\IEEEoverridecommandlockouts                

\usepackage{acronym}
\usepackage{pgfgantt}
\usepackage[utf8]{inputenc}
\usepackage{amsmath}
\usepackage{amssymb}
\usepackage{bm}
\usepackage[ruled,vlined]{algorithm2e}
\usepackage[final]{pdfpages}

\usepackage{amsthm}
\usepackage{dirtytalk}
\usepackage{mathtools}
\usepackage{tikz}
\usetikzlibrary{arrows,automata}
\usepackage{romannum}
\usepackage{wrapfig}
\usepackage{todonotes}
\usepackage{siunitx}
\usepackage{lineno}
\usepackage{bbm}
\usepackage{subcaption}
\usepackage{tikz}
\usepackage{float}
\usepackage{graphicx}
\usepackage{color}
\usepackage{url}
\usepackage{multirow}
\usepackage{comment}
\usepackage{paralist}
\usepackage{scrextend}
\usepackage{amsfonts}

\newcommand{\ie}{{\it i.e.}}

\newcommand{\reals}{\mathbb{R}}
\newcommand{\naturals}{\mathbb{N}}

\newcommand{\Expect}{\mathop{\bf E{}}}

\newcommand{\norm}[1]{\lVert #1 \rVert}
\newcommand{\abs}[1]{\lvert#1\rvert}

\newcommand{\BoolTrue}{\mbox{\sf true}}

\newcommand{\Always}{\Box}
\newcommand{\Eventually}{\Diamond}
\newcommand{\Next}{\bigcirc }

\newcommand{\Until}{\mbox{$\, {\sf U}\,$}}

\newcommand{\indicator}{\mathbf{1}}
\newcommand{\goal}{\mathsf{goal}}

\newcommand{\calAP}{\mathcal{AP}}
\newcommand{\calL}{\mathcal{L}}
\newcommand{\calA}{\mathcal{A}}

\newcommand{\calM}{\mathcal{M}}
\newcommand{\calX}{\mathcal{X}}

\newcommand{\calK}{\mathcal{K}}
\newcommand{\calB}{\mathcal{B}}
\newcommand{\nat}{\mathbb{Z}}

\newcommand{\inv}{\mathsf{\bf Inv}}
\newcommand{\guard}{\mathsf{\bf Guard}}

\newcommand{\level}{\mathcal{L}}

\acrodef{pa}[PA]{Probabilistic Automaton}
\acrodef{dfa}[DFA]{Deterministic Finite Automaton}
\acrodef{nfa}[NFA]{Nondeterministic Finite Automaton}
\acrodef{sccs}[SCCs]{Strongly Connected Components}
\acrodef{mdp}[MDP]{Markov Decision Process}
\acrodef{mdps}[MDPs]{Markov Decision Processes}
\acrodef{rl}[RL]{Reinforcement Learning}
\acrodef{tvi}[TVI]{Topological Value Iteration}
\acrodef{vi}[VI]{Value Iteration}
\acrodef{cpu}[CPU]{Central Processing Unit}
\acrodef{ltl}[LTL]{Linear Temporal Logic}
\acrodef{ap}[AP]{atomic propositions}
\acrodef{adp}[ADP]{Approximate Dynamic Programming}
\acrodef{tadp}[TADP]{Topological Approximate Dynamic Programming}
\acrodef{scltl}[sc-LTL]{syntactically co-safe LTL}
\acrodef{vi}[VI]{Value Iteration}
\acrodef{kl}[KL]{Kullback–Leibler}

\theoremstyle{definition}
\newtheorem{definition}{Definition}[section]
\newtheorem{theorem}{Theorem}[section]

\newtheorem{lemma}[theorem]{Lemma}

\newtheorem{remark}{Remark}
\newtheorem{problem}{Problem}
\newtheorem{proposition}{Proposition}[section]
\newtheorem{example}{Example}[section]

\title{\LARGE \bf
Topological Approximate Dynamic Programming under Temporal Logic Constraints
}

\author{Lening Li and Jie Fu
\thanks{L. Li and J. Fu are with Robotics Engineering Program, Department of Electrical and Computer Engineering,
Worcester Polytechnic Institute, Worcester, MA, 01609, USA,
{\tt\small lli4, jfu2@wpi.edu}}%
}

\begin{document}

\maketitle
\thispagestyle{empty}
\pagestyle{empty}

\begin{abstract}
In this paper, we develop a model-free approximate dynamic programming method for stochastic systems modeled as Markov decision processes to maximize the probability of satisfying high-level system specifications expressed in a subclass of temporal logic formulas---syntactically co-safe linear temporal logic. Our proposed method includes two steps: First, we decompose the planning problem into a sequence of sub-problems based on the topological property of the task automaton which is translated from a  temporal logic formula. Second, we extend a model-free approximate dynamic programming method to solve value functions, one for each state in the task automaton, in an order reverse to the causal dependency. Particularly, we show that the run-time of the proposed algorithm does not grow exponentially with the size of specifications. The correctness and efficiency of the algorithm are demonstrated using a robotic motion planning example.
\end{abstract}

\section{INTRODUCTION}
Temporal logic is a formal language to describe desired system properties, such as safety, reachability, obligation, stability, and liveness \cite{manna2012temporal}. 
This paper introduces a model-free \ac{rl} method for stochastic systems modeled as \ac{mdps}, where the planning objective is to maximize the (discounted) probability of satisfying constraints expressed in a subclass of temporal logic---\ac{scltl} formulas~\cite{belta2017temporal}.

Various model checking and probabilistic verification methods for \ac{mdp} are model-based~\cite{ding2014optimal,baier2008principles}. For systems without a model but with a blackbox simulator, \ac{rl} methods for \ac{ltl} constraints have been developed with both model-based and model-free methods \cite{fu2014probably,wen2016probably, alshiekh2018safe}. A model-based \ac{rl} learns a model and a near-optimal policy simultaneously. A model-free \ac{rl} learns only the near-optimal policy from sampled trajectories in the stochastic systems. However, model-free \ac{rl} methods, such as policy gradient and actor-critic methods~\cite{sutton1999policy,konda1999actor,sutton2018reinforcement}, face challenges when being used for planning with temporal logic constraints: a \ac{ltl} formula is translated into a sparse reward signal. The learner receives a reward of $1$ if the constraint is satisfied. This sparse reward provides little gradient information in the policy/value function search. The problem becomes more severe when complex specifications are involved. Consider the following example, a robot needs to visit regions A, B, and C, but if it visits D, then it must visit B before C. If the robot only visits A or B, then it will not receive any reward. When the state space of the \ac{mdp} is large, a learner not receiving any reward has no way to improve its current policy. To address reward sparsity, reward shaping~\cite{ng1999policy} has been developed. Reward shaping introduces additional reward signals while guaranteeing the policy invariance---the optimal policy remains the same with/without shaping. However, this method has strict requirements for the range of shaping potentials, which is hard to define when \ac{ltl} constraints are considered. 

In this work, we propose a different approach to mitigate the challenges in \ac{rl} under sparse reward signals in \ac{ltl} than reward shaping. Our approach is inspired by an idea for efficient value iteration: In an acyclic \ac{mdp}, there exists an optimal backup order, such that each state in the \ac{mdp} only needs to perform one-step backup operation in value iteration~\cite{bertsekas1995dynamic}. In \cite{dai2011topological}, the authors generalize this optimal backup order for acyclic \ac{mdp}s to general \ac{mdp}s. They develop a \ac{tvi} method that divides an \ac{mdp} into \ac{sccs} and then solves the values of states for each component sequentially in the topological order. Although it seems straightforward to apply \ac{tvi} to the product \ac{mdp}, which is obtained by augmenting the original \ac{mdp} with a finite set of memory states related to the task, the solution suffers from scalability issue. This may be mitigated with the use of \ac{adp}~\cite{li2018approximate}. The \ac{adp} makes the use of value function approximations to approximate the optimal solutions of large-scale problems. Hence, we propose a \emph{Topological Approximate Dynamic Programming} (TADP) method that includes two stages: Firstly, we translate the task formula into a \ac{dfa} referred as the \emph{task \ac{dfa}}, and then exploit the graphical structure in the automaton to determine a topological optimal backup order for \emph{a set of value functions}---one for each discrete state in the task \ac{dfa}. Value functions are related by the transitions in the task \ac{dfa} and jointly determine the optimal policy based on the Bellman equation. Secondly, we introduce function approximations for the set of value functions to reduce the number $N$ of decision variables---the number of states in the product \ac{mdp}---to a number $M$ of weights in function approximations, where $M\ll N$. Finally, we integrate a model-free \ac{adp} with the backup ordering to solve the set of value function approximations, one for each task state, in an optimal order. By doing this, the sparse reward received when the task is completed is propagated back to earlier stages of task completion, which provides meaningful gradient information for the learning algorithm. 

Exploiting the structure of task \ac{dfa}s for planning has been considered in~\cite{schillinger2018auctioning} where the authors partition the task \ac{dfa} into \ac{sccs} and then define progress levels towards satisfaction of the specification. In this work, we formally define a topological backup order based on the causal dependency among states in a task \ac{dfa}. We prove the optimality in this backup order. Further, this backup order can be integrated with the actor-critic method for \ac{ltl}-constrained planning in \cite{wang2015temporal} or other \ac{adp} methods that solve value function approximations to address the sparse reward problem.

The rest of the paper is structured as follows. Section~\ref{sec:preliminaries} provides some preliminaries. Section~\ref{sec:main_result} contains the main results of the paper, including computing the topological order, proof of optimality in this order, and the \ac{tadp} algorithm. The correctness and effectiveness of the proposed method are experimentally validated in the Section~\ref{sec:case_study} with a robotic motion planning example. Section~\ref{sec:conclusion} summarizes.

\section{PRELIMINARIES}\label{sec:preliminaries}

Notation: Given a finite set $X$, let $\Delta(X)$ be the set of probability distributions over $X$. The size of the set $X$ is denoted by $\abs{X}$. Let $\Sigma$ be an alphabet (a finite set of symbols). Given $k \in \nat^{+}$, $\Sigma^k$ indicates a set of words with length $k$, $\Sigma^{\leq k}$ indicates a set of words with length smaller than or equal to $k$, and $\Sigma^0 = \lambda$ is the empty word. $\Sigma^\ast$ is the set of finite words (also known as Kleene closure of $\Sigma$), and $\Sigma^\omega$ is the set of infinite words.  $\indicator_X$ is the indicator function with $\indicator_X(x) =1, \text{ if } x \in X$ and $0$ otherwise.

\subsection{Syntactically co-safe Linear Temporal Logic}
Syntactically co-safe \ac{ltl} formulas~\cite{Kupferman2001} are a well-defined subclass of \ac{ltl} formulas. Given a set of atomic propositions $\calAP$, the syntax of \ac{scltl} formulas is defined as follows:
\begin{equation*}
    \varphi \coloneqq \BoolTrue \mid p \mid \neg p \mid \varphi_1 \wedge \varphi_2 \mid \Next \varphi \mid \varphi_1 \Until \varphi_2,
\end{equation*}
where $\varphi, \varphi_1$ and $\varphi_2$ are \ac{scltl} formulas, $\BoolTrue$ is the unconditional true, and $p$ is an atomic proposition. Negation ($\neg$), conjunction ($\wedge$), and disconjunction ($\vee$) are standard Boolean operators. A \ac{scltl} formula can contain temporal operators \say{Next} ($\Next$), \say{Until} ($\Until$), and \say{Eventually} ($\Eventually$). However, temporal operator \say{Always} ($\Always$) is not contained in \ac{scltl}.

An infinite word with alphabet $2^\calAP$ satisfying a \ac{scltl} formula always has a finite-length good \emph{prefix}~\cite{Kupferman2001}~\footnote{$u$ is a prefix of $w$, \ie, $w=uv$ for $u, v\in\Sigma^\ast$}. Formally, given a \ac{scltl} formula $\varphi$ and an infinite word $w = r_0 r_1 \cdots $ over alphabet $2^{\calAP}$, $w \models \varphi$ if there exists $n \in \naturals$, $w_{[0:n]} \models \varphi$, where $w_{[0:n]} = r_0 r_1 \cdots r_n$ is the length $n+1$ prefix of $w$. Thus, an \ac{scltl} formula $\varphi$ over $2^\calAP $ can be translated into a \ac{dfa} $\calA_{\varphi} = \langle Q, \Sigma, \delta, q_0, F \rangle$, where $Q$ is a finite set of states,  $\Sigma = 2^{\calAP}$ is a finite set of symbols called the alphabet, $\delta: Q \times \Sigma \rightarrow Q$ is a transition function, $q_0\in Q$ is  an initial  state, and $F\subseteq Q$ is a set of accepting states. A transition function is recursively extended in the general way: $\delta(q,aw)=\delta( \delta(q,a),w )$ for given $a\in \Sigma$ and $w\in \Sigma^\ast$. A word $w = uv$ is \emph{accepting} if and only if $\delta(q, u)\in F$. \ac{dfa} $\calA_\varphi$ accepts the set of words satisfying $\varphi$.

We consider stochastic systems modeled as \ac{mdps}. We introduce a labeling function to relate paths in an \ac{mdp} $M$ to a given specification described by an \ac{scltl} formula.

\begin{definition}[Labeled \ac{mdp}]
\label{def:labeled_mdp}
	A labeled \ac{mdp} is a tuple $M = \langle S, A, s_0, P, \calAP, L \rangle$, where $S$ and $A$ are finite state and action sets, $s_0$ is the initial state, the transition probability function $P(\cdot \mid s, a) \in \Delta(S)$ is defined as a probability distribution over the next state given action is taken at the current state, $\calAP$ denotes a finite set of atomic propositions, and $L: S \rightarrow 2^{\calAP}$ is a labeling function mapping each state to the set of atomic propositions true in that state.
\end{definition}

A finite-memory stochastic policy in the \ac{mdp} is a function $\pi: S^{\ast} \rightarrow \Delta(A)$ that maps a history of state sequence into a distribution over  actions. A Markovian stochastic policy in the \ac{mdp} is a function $\pi: S \rightarrow \Delta(A)$ that maps the current state into a distribution over actions. Given an \ac{mdp} $M$ and a policy $\pi$, the policy induces a Markov chain $M^{\pi} = \{s_t \mid t = 0, \dots, \infty \}$ where $s_i$ is the random variable for the $i$-th state in the Markov chain $M^{\pi}$, and it holds that $s_{i+1} \sim P(\cdot \mid s_i, a_i)$ and $a_i \sim \pi(\cdot \mid s_0 s_1 \ldots s_i)$.

Given a finite (resp. infinite) path $\rho  =s_0s_1\ldots s_N \in S^\ast$ (resp. $\rho \in S^\omega$), we obtain a sequence of labels $L(\rho) = L(s_0)L(s_1)\ldots L(s_N)\in \Sigma^\ast$ (resp. $L(\rho) \in \Sigma^\omega$). A path $\rho$ satisfies the formula $\varphi$, denoted by $\rho \models \varphi$, if and only if $L(\rho) $ is accepted by $\calA_\varphi$. Given a Markov chain induced by policy $\pi$, the probability of satisfying the specification, denoted by $P(M^\pi\models \varphi)$, is the expected sum of the probabilities of paths satisfying the specification.
\[
P(M^\pi \models \varphi) \coloneqq \Expect\left[\sum_{t=0}^\infty \indicator(\rho_t\models \varphi)\right],
\]
where ${\rho}_t =s_0s_1\ldots s_t$ is a path of length $t+1$ in $M^\pi$. 


\begin{problem}\label{prob:maxprob}
Given a labeled \ac{mdp} $M$ and an \ac{scltl} formula $\varphi$, we can have a product \ac{mdp} $\calM$. The \emph{MaxProb} problem is to synthesize a policy $\pi$ that maximizes the probability of satisfying $\varphi$. Formally,
\begin{equation*}
\pi^\ast  =\arg\max_{\pi}P(M^\pi \models \varphi).
\end{equation*}
\end{problem}

\begin{definition}[Product \ac{mdp}]
Given a labeled \ac{mdp} $M = \langle S,A,s_0,P,\calAP, L\rangle$, an \ac{scltl} formula $\varphi$ with the corresponding \ac{dfa} $\calA_\varphi = \langle Q, \Sigma, \delta, q_0, F \rangle$, the product of $M$ and $\calA_{\varphi}$ is denoted by $M \otimes \calA_{\varphi} = \langle S \times Q ,(s_0, \delta(q_0,L(s_0)), S \times F, A, \delta, R \rangle$ with (1) a set of states, $S\times Q$, (2) an initial state, $(s_0, \delta(q_0,L(s_0))$, (3) the set of accepting states, $S \times F$, (4) the transition function defined by $ P(((s',q'),a') \mid (s, q), a) = P(s' \mid s, a) \indicator_{\{q'\}}(\delta(q,L(s)))$, (5) the reward function $R:S\times Q\times A\rightarrow [0,1]$. Formally,
\begin{align}
R((s,q),a) = \sum_{(s', q')}&P((s', q') \mid (s, q), a)\cdot \indicator_{F}(q').
\end{align}
We let all states in $ S\times F$  sink/absorbing states, \ie, for any $(s,q)\in S\times F$, for any $a\in A$, $P((s,q)|(s,q),a)=1$ and $R((s,q),a) =0$. For clarity, we denote this product \ac{mdp} by $\calM_{\varphi}$, \ie, $\calM_{\varphi} = M \otimes \calA_{\varphi}$. When the specification $\varphi$ is clear from the context, we denote the product \ac{mdp} by $\calM$.
\end{definition}

By definition, the path will receive a reward of $1$ if it ends in the set of accepting states $S \times F$. The total expected reward given a policy $\pi$ is the probability of satisfying the formula $\varphi$. By maximizing the total reward we find an optimal policy for the \emph{MaxProb} problem. In practice, we are often interested in maximizing a discounted total reward, which is the discounted probability of satisfying $\varphi$.

The planning problem is to solve the optimal value function and policy function satisfying
\begin{equation}
\label{eq:policy_softmax}
\begin{split}
& V((s, q))= \tau \log \sum_{a}\exp (R((s,q),a) \\
& \qquad + \gamma \sum_{s',q'}P((s',q') \mid (s,q), a)V((s',q'))) / \tau), \\
    & Q((s, q), a) = R((s,q),a) + \gamma \Expect_{(s',q')} V((s',q')),\\
    &\pi(a \mid (s, q)) = \exp((Q((s, q), a)- V((s, q)))/\tau),
    \end{split}
\end{equation}
where $\tau$ is an user-specified temperature parameter and $\gamma$ is a discounting factor. We use the softmax Bellman operator~\cite{nachum2017bridging} instead of the hardmax Bellman operator~\cite{sutton1998introduction}. \ac{vi} can solve the optimal value function in the product \ac{mdp} and converges in the polynomial time of the size of the state space, \ie, $\abs{S \times Q}$. However, \ac{vi} is model-based and difficult to scale to large planning problems with complex specifications. 

\section{MAIN RESULT}\label{sec:main_result}
We are interested in developing \emph{model-free} \ac{rl} algorithms for solving the \emph{MaxProb} problem. However, if we directly solve for approximately optimal policies in the product \ac{mdp} using the method in Section~\ref{sec:adp}, as the reward is sparse, it becomes a rare event to sample a path satisfying the specification. As a consequence, the estimate of the gradient in~\cite{li2018approximate} has a high variance with finite samples. To address this problem, we develop Topological Approximate Dynamic Programming (TADP) that leverages the structure property in the task automaton to improve the convergence due to sparse and temporally extended rewards with \ac{ltl} specifications.

\subsection{Hierarchical decomposition and causal dependency}
First, it is observed that given temporally extended goals, we can partition the product state space based on the discrete automaton states referred as discrete modes. The following definitions are generalized from almost-sure invariant set \cite{froyland2005statistically} in Markov chains to that in \ac{mdp}s. 
\begin{definition}[Invariant set and guard set]\label{def:inv_set_guard_set}
Given a \ac{dfa} mode $q\in Q$  and an \ac{mdp} $M$, the invariant set of $q$ with respect to $M$, denoted by $\inv(q, M)$, is a set of \ac{mdp} states such that no matter which action is selected, the system has probability one to stay within the mode $q$. Formally, 
\begin{align}
	\inv(q, M) = & \{s\in S\mid \forall a\in A,\forall s'\in S, P(s' \mid s,a) >0 \nonumber \\ & \implies \delta(q, L(s'))=q \},     
\end{align}
where $\implies$ means implication.

Given a pair $(q,q')$ of \ac{dfa} states, where $q \neq q'$, the set of \emph{guard states} of the transition from $q$ to $q'$, denoted by $\guard(q,q', M)$, is a subset of $S$ in which a transition from $q$ to $q'$ may occur. Formally,  
\begin{align}
	&\guard(q,q',M) =  \{s \in S \mid \exists a\in A, \exists s' \in S, \nonumber \\
	                 & P(s' \mid s,a) >0 \land \delta(q,L(s'))=q'\}. 
\end{align}
\end{definition}
When the \ac{mdp} $M$ is clear from the context, we refer $\inv(q,M)$ to $\inv(q)$ and $\guard(q,q',M)$ to $\guard(q,q')$.   

Next, we define \emph{causal dependency} between modes: In the product \ac{mdp} $\calM$, a state $(s_1,q_1)$ is \emph{causally dependent} on state $(s_2,q_2)$, denoted by $(s_1,q_1) \rightarrow (s_2,q_2)$, if there exists an action $a \in A$ such that
$P((s_2,q_2) \mid (s_1,q_1),a)>0$. This causal dependency is initially introduced in \cite{dai2011topological} and generalized to the state space of the product \ac{mdp}. 

According to Bellman equation \eqref{eq:bellman}, if there exists a probabilistic transition from $(s_1,q_1)$ to $(s_2,q_2)$ in the product \ac{mdp}, then the value $V(s_1,q_1)$ depends on the value $V(s_2,q_2)$.   

Two states can be causally dependent on each other. In that case, we say that these two states are \emph{mutually causally dependent.} Next, we lift the causal dependency from product \ac{mdp} to the task \ac{dfa}, by introducing casually dependent modes.

\begin{definition}[Causally dependent modes]\label{def:causally_dependent}
	A mode $q_1$ is \emph{causally dependent} on mode $q_2$ if and only if $\guard(q_1,q_2) \ne \emptyset$, that is, there exists a transition in the product \ac{mdp} from a state in mode $q_1$ to a state in mode $q_2$. A pair of modes $(q_1,q_2)$ is mutually causally dependent if and only if $q_1$ is causally dependent on $q_2$ and $q_2$ is causally dependent on $q_1$.
\end{definition}

\begin{definition}[Meta-mode]
	A \emph{meta mode} $X \subseteq Q$ is a subset of modes that are mutually causally dependent on each other. If a mode $q$ is not mutually causally dependent on any other modes, then the set $\{q\}$ itself is a meta mode. A meta mode $X$ is \emph{maximal} if there is no other mode in $Q\setminus X$ that is mutually causally dependent on a mode in $X$. 
\end{definition}

\begin{definition}[The maximal set of Meta-modes]
	$\cal X$ is said to be the set of \emph{maximal}  meta modes in the product \ac{mdp} if and only if it satisfies: i) any set $X \in \cal X$ is a maximal meta mode, ii) the union of sets in $\cal X$ yields the set $Q$, \ie,  $\cup_{X\in \mathcal{X}}X=Q$.
\end{definition}

\begin{lemma}
	The maximal set $\cal X$ of meta modes is a partition of $Q$.
\end{lemma}
\begin{proof}
	By the way of contradiction, suppose $\cal X$ is not a partition of $Q$, then there exists a mode $q\in X\cap X'$. Because $q$ is mutually causally dependent on all modes in $X$ as well as $X'$, then any pair $(q_1,q_2)\in X\times X'$ will be mutually causally dependent---a contradiction to the definition of $\cal X$.
\end{proof}

We denote $X \rightarrow X'$ if a mode $q \in X$  is causally dependent on mode $q' \in X'$. By the transitivity property, if $X_1\rightarrow X_2$ and $X_2\rightarrow X_3$, then we denote the causal dependency of $X_1$ on $X_3$ by $X_1\rightarrow^+ X_3$. The following lemma states that two states in the product \ac{mdp} are causally dependent if their discrete modes are causally dependent.

\begin{lemma}
\label{lma2}
    Given two meta-modes $X, X' \in \calX$, if $X \rightarrow^+ X'$ but not $X' \rightarrow^+ X$, then for any state $(s,q) \in S\times X$ and $(s',q')\in S\times X'$, it is the case that either $(s,q)\rightarrow^+(s',q')$ or these two states are causally independent.
\end{lemma}
\begin{proof}
    By the way of contradiction, if $(s,q)$ and $(s',q')$ are causally dependent and $(s',q')\rightarrow^+ (s,q)$, then there must exist a state $(s'', q'')$ such that $(s',q')\rightarrow^+ (s'',q'')$ and $(s'',q'')\rightarrow (s,q)$. Relating the causal dependency of states in the product \ac{mdp} and the definition of guard set, we have $s'' \in \guard(q'',q)$ and $q''\rightarrow q$. Further $q'\rightarrow^+ q''\rightarrow q$, thus $X'\rightarrow^+ X$, which is a contradiction as $X' \not \rightarrow^+ X$.
\end{proof}

Lemma~\ref{lma2} provides structural information about topological value iteration (TVI). If $X\rightarrow^+ X'$ and $X' \not \rightarrow^+ X$, then based on \ac{tvi}, we shall update the values for states in the set $\{(s,q)\mid q\in X'\}$ before updating the values for states in the set $\{(s,q)\mid q\in X\}$.

However, the causal dependency in meta-modes does not provide us with a total order over the set of maximally meta modes because two meta modes can be causally independent. A total order is needed for deciding the order in which the optimal value functions for modes are updated.
 
To obtain a total order, we construct a total ordered sequence of sets of maximal meta modes. Given the set of maximal meta-modes $\cal X$, 
\begin{enumerate}
    \item Let $\level_0 = \{ X \in {\cal X} \mid X\cap F\ne \emptyset\}$ and $i=1$.
    \item Let $\level_i= \{X\in {\cal X} \setminus \cup_{k=0}^{i-1} \level_k \mid \exists X'\in \level_{i-1}\text{ such that } X\rightarrow X'\}$, and increase $i$ by 1
    \item Repeat until $i=n$ and $\level_{n+1} = \emptyset$.
\end{enumerate}
We refer $\{\level_i \mid i=0,\ldots, n\}$ as \emph{level sets over meta modes}. Based on the definition of set $\{\level_i \mid i=0,\ldots, n\}$, we let $\rightsquigarrow$ define an ordering on level sets as follows: $\level_{i}\rightsquigarrow \level_{i-1} \mid i = 1, \ldots, n$.
We give the following two statements.

\begin{lemma}
    If there exists $X \in \calX$ such that $X \not\in \level_i$ for any $i=0,\ldots, n$, then the set of states in $X$ is not coaccessible from the final set $F$ of states in \ac{dfa} $\calA_\varphi$.
\end{lemma}
\begin{proof}
By construction, this meta mode $X$ is not causally dependent on any meta mode that contains $F$. Thus, it is not coaccessible in the task \ac{dfa} $\calA_\varphi$, \ie, there does not exist a word $w$ such that $\delta(q,w)\in F$ for some $q\in X$.
\end{proof}
If a \ac{dfa} is coaccessible, then we have $\cup_{i=0}^n \level_i  = \calX$. A state $q$ that is not coaccessiable from the final set $F$ should be trimmed  before planning because the value $V(s,q)$ for any $s\in S$ will not be used for optimal planning in the product \ac{mdp} to reach $F$. 

\begin{proposition}
    The ordering $\rightsquigarrow $ is a total order:
    \[
    \level_{n}\rightsquigarrow \level_{n-1}\ldots \rightsquigarrow \level_0.
    \]
\end{proposition}

\begin{theorem}[Optimal Backup Order~\cite{bertsekas1995dynamic}]
	If an \ac{mdp} is acyclic, then there exists an optimal backup order. By applying the optimal order, the optimal value function can be found with each state needing only one backup.
\end{theorem}

We generalize the Optimal Backup Order on an acyclic \ac{mdp} to the product \ac{mdp} as the following:

\begin{theorem}[Generalized optimal backup order for hierarchical planning]
\label{thm:backuporder} Given the  optimal planning problem in the product \ac{mdp} and the causal ordering of meta modes, by updating the value functions of all meta-modes in the same level set, in a sequence reverse to the ordering $\rightsquigarrow$, the optimal value function for each meta mode can be found with only one backup, \ie, solving the value function of that meta mode using value iteration or an \ac{adp} method that solves the value function approximation.
\end{theorem}
\begin{proof}We show this by induction. Suppose there exists only one level set, the problem is reduced to optimal planning in a product \ac{mdp} with only one update for value functions of meta-modes in this level set. When there are multiple level sets, each time the optimal planning performs value function update for one level set. The value $V(s,q)$ for $q\in X$  only depends on the values of its descent states, that is, $\{V(s',q') \mid (s,q)\rightarrow (s',q')\}$. It is noted that the mode $q'$ of any descendant $(s,q)$ must belong to either meta-modes $X$, or some $X' \in \calX$ such that $X \rightarrow^+ X'$. By definition of level sets, if $X\in \calL_i$, then $X' \in \calL_k$ for some $k \le i$. It means the value $V(s',q') $ for any descendant $(s,q)$ is either updated in level $\calL_k$, $k<i$, or along with the value $V(s,q)$, when $k=i$. As a result, after the value functions $\{V(\cdot, q)\mid q\in X, X\in\calL_i\}$ converge, the Bellman residuals of states in $\{(s,q)\mid q \in X, X\in \calL_{k}, k \le i\}$ remain unchanged, while the value functions of meta-modes in other level sets with higher levels are updated. Thus, each mode only needs to be updated once.
\end{proof}
\begin{example} \label{ex:specificition}
We use a simple example to illustrate. Given a system-level specification: $\Eventually(b \land \Next \Eventually c) \land \Eventually(a \land \Next \Eventually d)$, the corresponding \ac{dfa} is shown in Fig.~\ref{fig:automaton}. In this \ac{dfa}, each state is its own meta mode $X_{i}=\{q_{i+1}\}, i = 0, \ldots 8$. Different level sets $\calL_{i}, i = 0, \ldots, 4$, are contained in different styled ellipses.
\begin{figure}[!htb]
\centering
\vspace{-3ex}
	\begin{tikzpicture}[->,>=stealth',shorten >=1pt,auto,node distance=3cm,scale=0.5,semithick, transform shape]
		\draw[ultra thick,dashed,orange]        (0,0)       circle [x radius=2cm, y radius=0.75cm,];
        \draw[ultra thick,densely dotted,black] (0,-2)      circle [x radius=3cm, y radius=1cm,];
        \draw[ultra thick,solid,cyan]           (0,-4.25)   circle [x radius=5cm, y radius=1cm,];
        \draw[ultra thick,dashdotted,red]       (0,-6.5)    circle [x radius=3cm, y radius=1cm,];
        \draw[ultra thick,loosely dotted,blue]  (0,-8.5)   circle [x radius=2cm, y radius=0.75cm,];

        \node[] at (2,0)   {\LARGE $\calL_4$};
        \node[] at (3.5,-2)     {\LARGE $\calL_3$};
        \node[] at (5.5,-4.25)    {\LARGE $\calL_2$};
        \node[] at (3.5,-6.25)     {\LARGE $\calL_1$};
        \node[] at (3,-8.25)     {\LARGE $\calL_0$};

		\tikzstyle{every state}=[fill=black!10!white]
		\node[initial,state] (1)                        {$q_1$};
        \node[state] (2) [below left of=1]              {$q_2$};
        \node[state] (3) [below right of=1]             {$q_3$};
        \node[state] (4) [below left of=2]              {$q_4$};
        \node[state] (5) [below right of=2]             {$q_5$};
        \node[state] (6) [below right of=3]             {$q_6$};
        \node[state] (7) [below right of=4]             {$q_7$};
        \node[state] (8) [below right of=5]             {$q_8$};
        \node[state,accepting] (9) [below right of=7]   {$q_9$};
		\path[->]   
		(1) edge node {$a$}                 (2)
        (1) edge node {$b$}                 (3)
        (1) edge[loop above] node {$\neg a \wedge \neg b$}  (1)
        (2) edge node {$d$}                 (4)
        (2) edge node {$b$}                 (5)
        (2) edge[loop above] node {$\neg b \wedge \neg d$}  (2)
        (3) edge node {$a$}                 (5)
        (3) edge node {$c$}                 (6)
        (3) edge[loop above] node {$\neg a \wedge \neg c$}  (3)
        (4) edge node {$b$}                 (7)
        (4) edge[loop above] node {$\neg b$}  (4)
        (5) edge node {$d$}                 (7)
        (5) edge[loop right] node {$\neg c \wedge \neg d$}  (5)
        (5) edge node {$c$}                 (8)
        (6) edge node {$a$}                 (8)
        (6) edge[loop above] node {$\neg a$}  (6)
        (7) edge node {$c$}                 (9)
        (7) edge[loop above] node {$\neg c$}  (7)
        (8) edge node {$d$}                 (9)
        (8) edge[loop above] node {$\neg d$}  (8)
        (9) edge[loop below] node {$\BoolTrue$}  (9)
		;
	\end{tikzpicture}
	\caption{\ac{dfa} with respect to $\Eventually(b \land \Next \Eventually c) \land \Eventually(a \land \Next \Eventually d)$. }
	\label{fig:automaton}
\vspace{-3ex}
\end{figure}
\end{example}

\subsection{Model-free  \ac{adp} for planning with temporal logic constraints}\label{sec:adp} 
\ac{adp} makes use of value function approximations to solve for Problem~\ref{prob:maxprob}. First, let's define the softmax Bellman operator by
\begin{multline}
	 \calB V(s,q) = \tau \log \sum_{a} \exp ( (R((s,q), a) \\
	       + \gamma \sum_{(s',q')}P((s',q') \mid (s,q), a)V(s',q') ) / \tau ), 
	\label{eq:bellman}
\end{multline}
where $\tau > 0$ is an user-specified temperature parameter.

We introduce a mode-dependent value function approximation as follows: For each $q\in Q$, the value function is approximated by $V(\cdot; \theta_{q}): S  \rightarrow \reals$ where $\theta_{q} \in \reals^{\ell_{q}}$ is a parameter vector of length $\ell_q$. We use a linear function approximation of $V(\cdot;\theta_q) =\sum_{k=1}^{\ell_q} \phi_{k,q}(\cdot) \theta_{q}[k] = \Phi_q\theta_{q}$, where $\phi_{k,q}: S \rightarrow \reals, k = 1, \dots, \ell_{q}$ are pre-selected basis functions. We first define two sets: For meta-modes $X, X'\in \calX$, let 
 \[
\inv(X) = \bigcup_{q\in X}\inv(q), \quad \text{ and}
\]
 \[
\guard(X,X')=\bigcup_{q\in X,q'\in X'} \guard(q,q').\]

Given the level sets $\{\level_i \mid i=0,\ldots, n\}$, the computation of value function approximation for each \ac{dfa} mode is carried out in the order of level sets.
\begin{enumerate}
    \item Starting with level $0$, let $i = 0$ and $V((s,q);\theta_q)= 1$ (see Remark 2) for all $s\in S$ and $q\in F$. For each $X\in \level_0$, we solve an \ac{adp} problem:
      \begin{align}
	\min_{\{\theta_{q} \mid q\in 
	X \setminus F \}} & \sum_{(s, q) \in S \times X} c(s, q)V((s, q);  \theta_{q}),                           \\
	\mbox{subject to: } & \calB V((s, q);\theta_{q}) - V((s, q); \theta_{q}) \leq 0 , \nonumber \\
	 & \forall s\in \inv(X)  \bigcup \left(  \cup_{X'\in\calX}\guard(X, X')\right),\nonumber
\end{align}
    where parameters $c(s,q)$ are state relevant weights.
    All states $\{(s,q)\mid q\in F\}$ are absorbing with values of $1$. 
    The reward function $R((s,q),a)=0$ for all $s\in S$, $q\in X$, and $a\in A$. After solving the set of value functions $\{V(s,q;\theta_q) \mid q\in X, X\in \level_0\}$. The solution  of this \ac{adp} is proven to be  a tight upper bound of the optimal value function \cite{li2018approximate}. See Appendix for more information about this \ac{adp} method.
    \item Let $i \leftarrow i+1$.
    \item At the $i$-th step, given the value $\{V(s,q; \theta_q)\mid q\in X\land X\in \level_k, k<i\}$, we solve, for each $X\in \level_i$, an \ac{adp} problem stated as follows:
    \begin{align}
\label{eq:topo_LTL}
	\min_{\{\theta_q\mid q\in X\}} & \sum_{(s, q) \in S \times X} c(s, q)V((s, q);  \theta_{q}),                           \\
	\mbox{subject to: } & \calB V((s, q);\theta_{q}) - V((s, q); \theta_{q}) \leq 0 , \nonumber \\
	 & \forall s\in  \inv(X ) \bigcup \left(  \cup_{X'\in\calX}\guard(X, X')\right),\nonumber 
\end{align}
where  $V((s',q');\theta_{q'})$ to be solved either has $q'\in X$ or $q'\in X'$ for some $X'\in \level_k$, $k < i$. Note that by Theorem~\ref{thm:backuporder}, the meta-mode $X'$, for which $\guard(X,X')$ is nonempty, cannot be in a level set higher than $i$.
When $q'\in X'$ and $X'\in \level_i$, then $\theta_{q'}$ is a decision variable for this \ac{adp}. When $q'\in X'$ and $X'\in \level_k$ for some $k <i$, then the value $V((s',q');\theta_{q'})$ is computed in previous iterations and substituted. A state $(s',q')$ whose value is determined in previous iterations is made absorbing in this iteration. The reward function $R((s,q),a)=0$ for all $s\in S$, $q\in X$, and $a\in A$.
\item Repeat step $2,3$ until $i=n$. Return the set $\{V(s,q;\theta_q)\mid q\in Q\}$. The  policy is computed using the softmax Bellman operator, defined in \eqref{eq:policy_softmax} by substituting the value function $V(s,q)$ with its approximation $V((s,q);\theta_q)$.
\end{enumerate}

\begin{remark}
The problem solved by \ac{adp} is essentially a stochastic shortest path problem~\cite{bertsekas1991analysis}. For such a problem, two approaches can be used: One is to fix the values of states to be reached and assign the reward to be zero. During value iteration, the value of the states to be reached will be propagated back to the values of other states. The aforementioned reward design and \ac{adp} formulation use the first approach.
Another way to introduce a reward function defined by $R((s, q),a) = \sum_{s', q'}P((s', q') \mid (s, q), a) R((s',q')) $, where $R(s',q')= V((s',q');\theta_{q'})$ if $q \in X$, $q' \in X'$, $X \rightarrow X'$
, and $R(s',q')=0$ otherwise.
\end{remark}

\begin{remark}
A value iteration using the softmax Bellman operator finds a policy that maximizes a weighted sum of total rewards and the entropy of policy (see~\cite{nachum2017bridging} for more details). When the value/reward is small, the total entropy of policies accumulated with the softmax Bellman operator overshadows the value given by the reward function. This is called the value diminishing problem. Thus, for both cases, when the value $V((s',q');\theta_{q'})$ of the state to be reached is small, we scale this value by a constant $\alpha$ to avoid the value diminishing problem. Given the nature of the \emph{MaxProb} problem, with a reward of 1 being assigned when the \ac{ltl} constraint is satisfied, we almost always need to amplify the reward to avoid the value diminishing problem.
\end{remark}

\section{CASE STUDY}\label{sec:case_study}
\begin{figure}[!htb]
\vspace{-2ex}
    \centering
    \includegraphics[width=0.8\columnwidth]{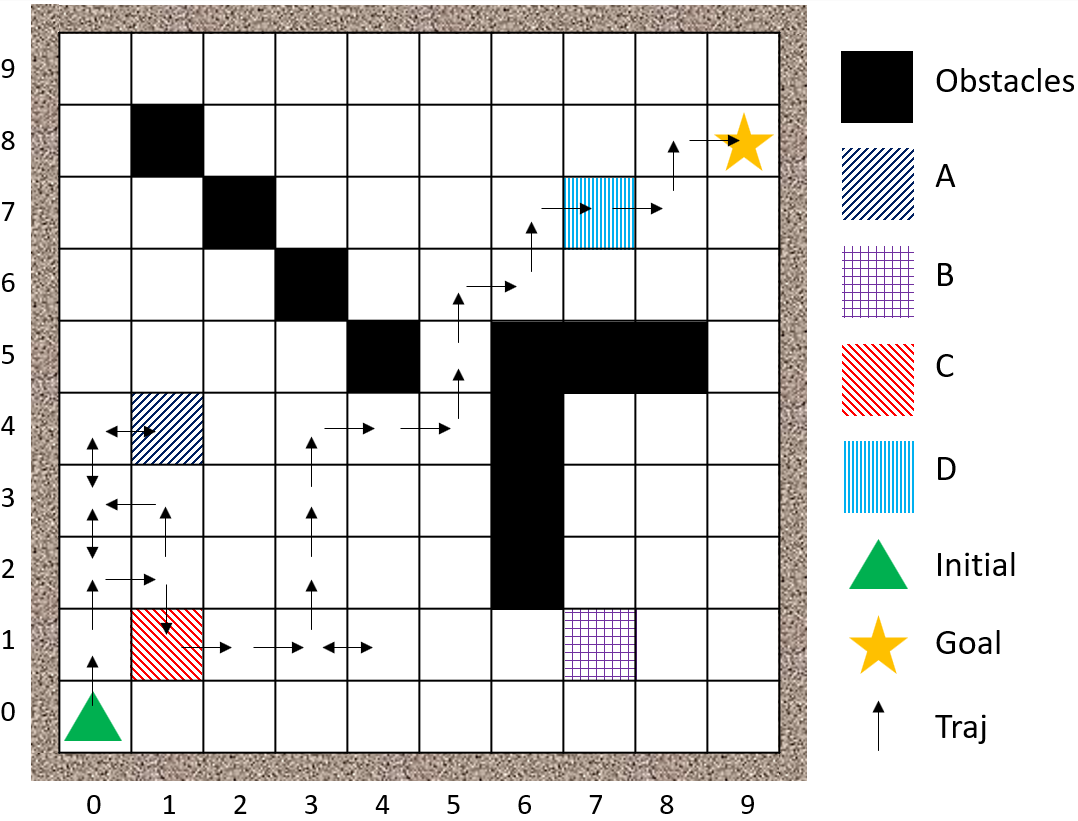}
    \caption{One simulation on the grid world.}
    \label{fig:ADP_sim}
\end{figure}

We validate the algorithm in a motion planing problem under a \ac{scltl} specification in a grid world. In this example, we consider the following specification: $\Eventually(((A \land (\neg B \Until C)) \lor (B \land (\neg A \Until D))) \land \Eventually \goal \land \Always \neg O)$, and the corresponding \ac{dfa} is plotted in Fig.~\ref{fig:case_automaton}. This specification describes that the system visits A, C, and goal sequentially, or the system visits B, D, and goal sequentially, while avoiding obstacles. Regions $A, B, C, D$, obstacles $O$, and goal are shown in Fig.~\ref{fig:ADP_sim}. The partitions of meta modes are shown in Fig.~\ref{fig:case_automaton} with different meta modes being boxed in different styled rectangles. The task automaton is partitioned into four meta modes $X_i, i = 0, \ldots, 3$, and each level set $\calL_{i}, i = 0, \ldots, 3$, contains one meta mode with the same index. The reward is defined as the following: the robot receives a reward of $60$ (an amplified reward to avoid value diminishing) if the trajectory satisfies the specification. In each state $s\in S$ and for robot's different actions (heading up ('U'), down ('D'), left('L'), right('R')), the probability of arriving at the \say{correct} cell is $1 - 0.03 \times \abs{N}$, and the probability of arriving a \say{wrong} cell is $0.03$, where $0.03$ is the randomness in the system and $\abs{N}$ is the number of possible succeeding states. We surround the gird world with walls. If the system hits the wall, it will be bounced back and stay at its original cell. All the obstacles are sink states, \ie, when a robot goes to an obstacle---it stays there forever.
\begin{figure}[!htb]
\vspace{-5ex}
\centering
	\begin{tikzpicture}[->,>=stealth',shorten >=1pt,auto,node distance=3cm,scale=0.5,semithick, transform shape]
	\draw[red,thick,dashed]             (-1,-1)     rectangle (0.8,1);
	\draw[blue,thick,densely dotted]    (1.2,-3)    rectangle(3,3);
	\draw[black,thick,solid]            (6.5,-1)    rectangle(8,1);
	
	\node[] at (8.5,-0.5)   {\LARGE $X_0$};
    \node[] at (5.5,-0.5)   {\LARGE $X_1$};
    \node[] at (3.5,-2.5)   {\LARGE $X_2$};
    \node[] at (0.0,-1.5)   {\LARGE $X_3$};
    
    \node[] at (8.5,1.5)    {\LARGE $\calL_0$};
    \node[] at (5,1.5)      {\LARGE $\calL_1$};
    \node[] at (3.5,3.0)    {\LARGE $\calL_2$};
    \node[] at (0.0,2.0)    {\LARGE $\calL_3$};
    
	\draw[green,thick,double]    (3.4,-1)  rectangle(5,1);
		\tikzstyle{every state}=[fill=black!10!white]
		\node[initial,state]    (0)                     {$q_1$};
        \node[state]            (1) [above right of=0]  {$q_2$};
        \node[state]            (2) [below right of=0]  {$q_3$};
        \node[state]            (3) [above right of=2]  {$q_4$};
        \node[state, accepting] (4) [right of=3]        {$q_5$};
		\path[->]   
		(0) edge                node            {$a$}                   (1)
        (0) edge                node[below]     {$b$}                   (2)
        (0) edge[loop above]    node            {$\neg b \wedge \neg a$}  (0)
        (1) edge[bend right]    node            {$b$}       (2)
        (1) edge                node            {$c$}       (3)
        (1) edge[loop above]    node            {$\neg b \wedge \neg c$}  (1)
        (2) edge[bend right]    node            {$a$}       (1)
        (2) edge                node[below]     {$d$}       (3)
        (2) edge[loop below]    node            {$\neg a \wedge \neg d$}  (2)
        (3) edge                node            {$\goal$}  (4)
        (4) edge[loop above]    node            {$\BoolTrue$}   (4)
		;
	\end{tikzpicture}
\caption{Automaton $\Eventually(((A \land (\neg B \Until C)) \lor (B \land (\neg A \Until D))) \land \Eventually \goal \land \Always \neg O)$. Meta modes and the level sets are marked.}
\label{fig:case_automaton}
\vspace{-2ex}
\end{figure}

The planning objective is to find an approximately optimal policy for satisfying the specification with a maximal probability. We compare the \ac{tadp} with \ac{vi} and \ac{tvi} to show the correctness and efficiency.

\paragraph*{Parameters}
We use the following parameters: the user-specified temperature is $\tau = 2$, discounting factor is $\gamma=0.9$, and error tolerance is $\epsilon=10^{-3}$. The tolerance is shared by \ac{tvi}, \ac{vi}, and \ac{tadp}, where the stopping criterion is $\|\max(V^{j}(s) - V^{j-1}(s))\| \leq \epsilon$ for $j$-th iteration of in \ac{tvi} and \ac{vi} and for each inner $j$-th iteration of \ac{tadp}, respectively.

We adopt the following initial parameters in the \ac{tadp} algorithm for the $k$-th problem (see appendix for the meanings of parameters): the coefficient of the penalty is $b=1.5$, the learning rate is $\eta = 0.1$, the penalty parameter is $\nu = 2.0$, and the Lagrangian multipliers is $\lambda = 0$. During each inner iteration, we sample $30$ trajectories of length $\leq 3$. The value function $V(\cdot; \theta_{q})$ is approximated by a weighted sum of Gaussian Kernels: $V(\cdot; \theta_q) = \Phi_{q} \theta_q$, where basis functions $\Phi_q= [\phi_1,\phi_2,\ldots, \phi_{\ell_q}]^\intercal$ are defined as the following: $\Phi_{j}(s)  = K(s, c^{(j)}) $ and $K(s, s') = \exp(-\frac{SP(s,s')^2}{2\sigma^2}), $ where $\{c^{(j)}, j=1,\ldots, \ell_q\} $ is a set of pre-selected centers and $\sigma=1$. In this example, we select the centers to be uniformly selected points with interval $1$ within the grid world.

After the \ac{tadp} converges, we obtain the policy from the converged value functions computed by \ac{tadp} and simulate the system. We plot one simulation of the system in Fig.~\ref{fig:ADP_sim}. The system starts at the initial state $s_{init}$, then it visits region $A$ then region $C$, eventually it visits the goal state $s_{goal}$.

\begin{figure*}[!htb]
\centering
    \begin{subfigure}[b]{0.25\textwidth}
    \includegraphics[width=\textwidth]{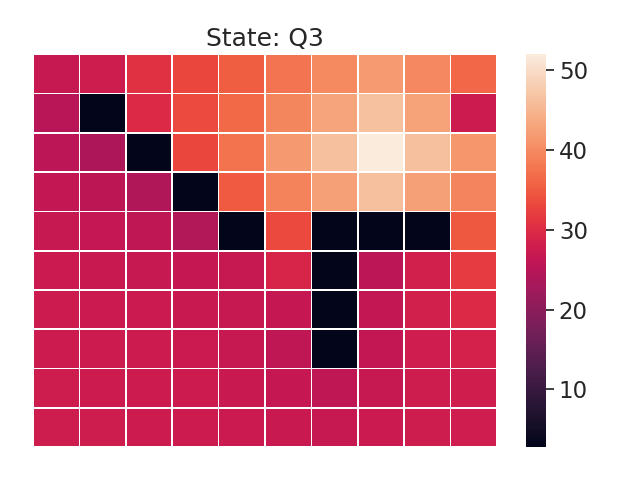}
    \caption{ }
    \label{fig:Product_Q3_heatmap}
    \end{subfigure}
    \begin{subfigure}[b]{0.25\textwidth}
    \includegraphics[width=\textwidth]{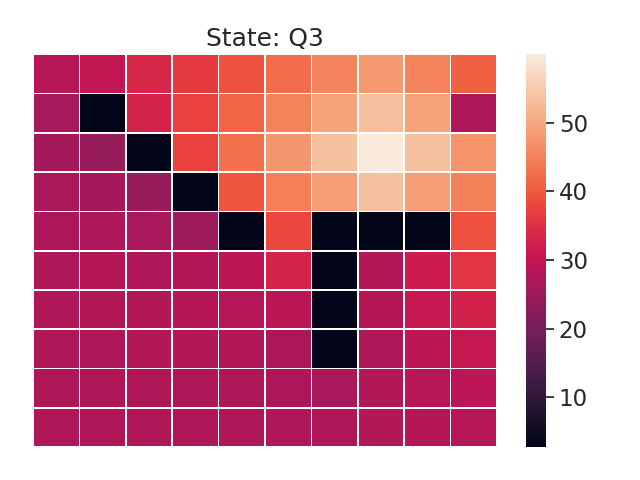}
    \caption{ }
    \label{fig:TVI_Q3_heatmap}
    \end{subfigure}
    \begin{subfigure}[b]{0.25\textwidth}
    \includegraphics[width=\textwidth]{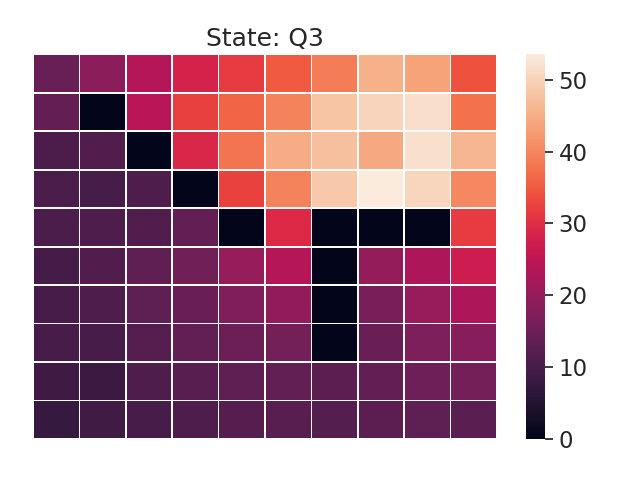}
    \caption{ }
    \label{fig:ADP_Q3_heatmap}
    \end{subfigure}
    \begin{subfigure}[!htb]{0.35\textwidth}    
    \includegraphics[width=0.3\textwidth]{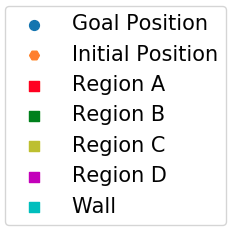}\hfill
    \includegraphics[width=0.7\textwidth]{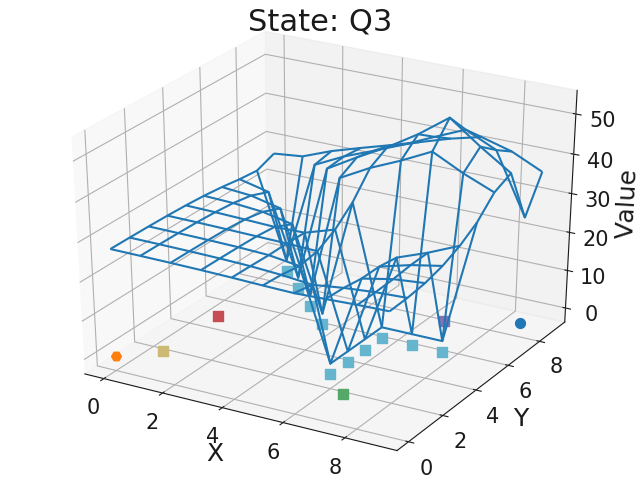}
    \caption{ }
    \label{fig:Product_Q3_value}
    \end{subfigure}
    \begin{subfigure}[!htb]{0.25\textwidth}
    \includegraphics[width=\textwidth]{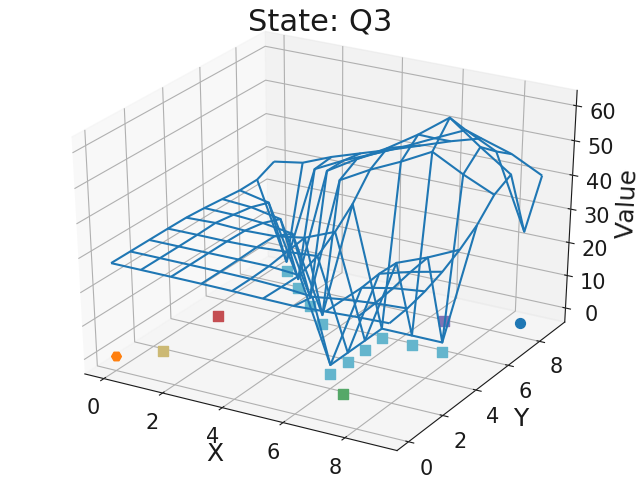}
    \caption{ }
    \label{fig:TVI_Q3_value}
    \end{subfigure}
    \begin{subfigure}[!htb]{0.25\textwidth}
    \includegraphics[width=\textwidth]{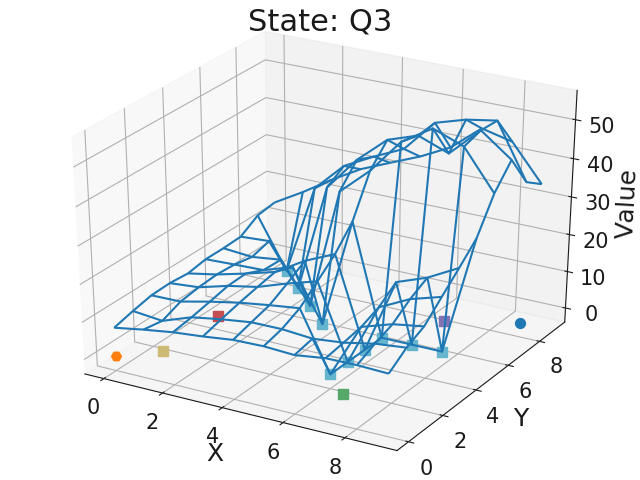}
    \caption{ }
    \label{fig:ADP_Q3_value}
    \end{subfigure}
    \caption{Comparison between \ac{vi}, \ac{tvi}, \ac{tadp} for different states at $q_3$: (a) (b) (c) are the heatmaps of $V(\cdot, q_3)$ obtained by \ac{vi}, \ac{tvi}, and \ac{tadp}, respectively.  (d) (e) (f) are the corresponding value surfs of $V(\cdot, q_3)$ obtained by \ac{vi}, \ac{tvi}, and \ac{tadp}, respectively.}
    \label{fig:comparison}
    \vspace{-4ex}
\end{figure*}

\paragraph*{Value Comparison}
In Fig.~\ref{fig:comparison}, we plot the heatmaps and values for different states at $q_3$ obtained by \ac{vi}, \ac{tvi}, and \ac{tadp}. In heatmaps, the brighter the area is, the higher value of that area is. The results from Fig.~\ref{fig:Product_Q3_heatmap} and ~\ref{fig:TVI_Q3_heatmap} shows that \ac{vi} and \ac{tvi} both show the most bright area is at $(7, 7)$. In Fig.~\ref{fig:ADP_Q3_heatmap}, the area around $(7,7)$ obtained has relatively bright color. The heatmap of \ac{tadp} is not exactly the same as the other two. This is due to the approximation error. Comparing three value surfs in Fig.~\ref{fig:Product_Q3_value}, ~\ref{fig:TVI_Q3_value}, and ~\ref{fig:ADP_Q3_value}, we are able to see the similarity between these three value surfs. 

\paragraph*{Run-time}
We conduct two experiments for different sizes of grid worlds, \ie, $10 \times 10$ and $20 \times 20$. We show the results in Table~\ref{table:running_time}. In different sizes of gird worlds, comparing \ac{vi} and \ac{tvi} run-times are reduced by $38.45\%$ and $53.62\%$ by exploiting the topological structure, and the total numbers of Bellman Backup Operations are reduced by $7.71\%$ and $7.76\%$. In terms of simple specifications, the decomposition occupies major CPU time, but
exploiting the topological structure will be leveraged if more complex specifications are associated. The \ac{tadp} converges after $135.96$ seconds and $1117.91$ seconds, respectively for different sizes of grid worlds. The run-times of the \ac{tvi} and \ac{vi} in a $20 \times 20$ grid world are $14$-$20$ times their run-times in the $10 \times 10$ grid world. However, the run-time of \ac{tadp} in a $20 \times 20$ grid world is only $8$ times the run-time of \ac{tadp} in the $10 \times 10$ grid world. 
\ac{tadp} is more beneficial in large \ac{mdp} problems or with more complex specifications. It is noted that though \ac{tadp} takes in general longer time to converge, but it is model-free. \ac{tvi} and \ac{vi} are model-based. 
\begin{table}[!htb]
\centering
\resizebox{\linewidth}{!}{%
\begin{tabular}{|c|c|c|c|c|}
\hline
\multicolumn{2}{|c|}{Algorithms} & VI & TVI & TADP \\ \hline
\multirow{2}{*}{$10 \times 10$} & Bellman Backup Operations (times) & 64620 & 59636 & N/A \\ \cline{2-5} 
& Run-time (Seconds) & 11.52 & 7.09 & 135.96 \\ \hline
\multirow{2}{*}{$20 \times 20$} & Bellman Backup Operations (times) & 280620 & 258836 & N/A \\ \cline{2-5} 
& Run-time (Seconds) & 222.56 & 103.21 & 1117.59 \\ \hline
\end{tabular}%
}
\caption{Bellman Backup Operations and Run-times between \ac{vi}, \ac{tvi}, and \ac{tadp}. Note that \ac{tvi} has a significantly shorter run-time.}
\label{table:running_time}
\vspace{-4ex}
\end{table}

\begin{figure}[!htb]
\centering
    \includegraphics[width=0.7\columnwidth]{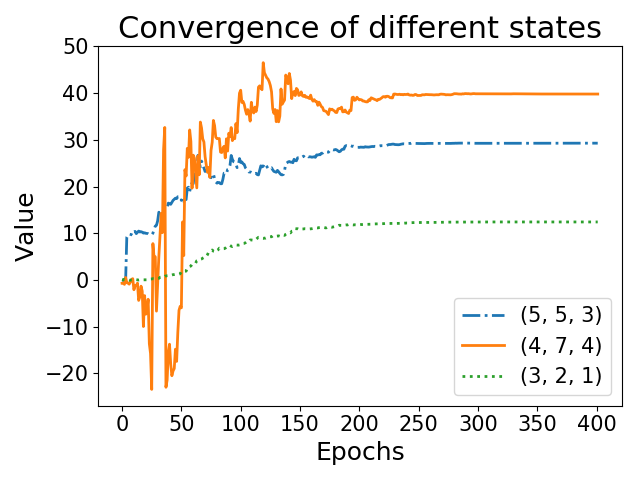}
    \caption{The convergence of values in \ac{tadp} in the $10 \times 10$ stochastic grid world for different states in the product \ac{mdp}. A product state $(5,5,3)$ means the grid cell $(5,5)$ and the \ac{dfa} state $q_3$.}
    \label{fig:convergence}
\vspace{-3ex}
\end{figure}

\paragraph*{Convergence}
In Fig.~\ref{fig:convergence}, we plot the convergence of values for different states in the $10\times 10$ grid and modes in automaton against epochs, which is the number of inner iterations in \ac{tadp}. It indicates that the values initially oscillate, but all values converge after $250$ iterations. 

\paragraph*{Statistical Results} We want to quantify and compare the performance between different methods. We update the policies from converged value functions computed by \ac{tadp} and \ac{tvi}. We simulate trajectories for $500$ times and compare the percentages of trajectories reaching the goal out of the total trajectories. We limit the max time-step to be $500$, that is, if the system cannot reach the goal within $500$ steps, then the system fails to reach the goal. Moreover, if the system reaches any sink states, then the system fails to reach the goal. Otherwise, the system successfully reaches the goal. We conduct two statistical experiments under the same setting with different starting potions, \ie, $(1,2,2)$ and $(2,2,4)$. Starting with $(1,2,2)$ the percentages of reaching the goal under \ac{tadp} and \ac{tvi} are $66.2\%$ and $86.8\%$. Starting with $(2,2,4)$ the percentages of reaching the goal under \ac{tadp} and \ac{tvi} are $80\%$ and $88.6\%$. The results indicate that the policy computed by \ac{tadp} is suboptimal due to the nature of \ac{adp}, but the performance gap between two policies is not significant. 

\section{CONCLUSION}\label{sec:conclusion}
We present a topological approximate dynamic programming method to maximize the probability of satisfying high-level system specifications in \ac{ltl}. We decompose the product \ac{mdp} and define the topological order for updating value functions at different task modes to mitigate the sparse reward problems in model-free \ac{rl} with \ac{ltl} objectives.  The correctness of the algorithm is demonstrated on a robotic motion planning problem under \ac{ltl} constraints. It is noted that one needs to update the value functions for all discrete states in a meta-mode at a time. When the size of meta-mode is large, then the number of parameters in value function approximations to be solved is large, which raises the scalability issue due to the complexity of the specifications. We will investigate action elimination technique within the framework, not at the low-level actions in the \ac{mdp}, but at the high-level decisions of transitions in the task \ac{dfa}. By eliminating transitions in the \ac{dfa}, it is possible to decompose large meta-mode into a subset of small meta-modes whose value functions can be efficiently solved.


\appendix
In this section, we present a model-free \ac{adp} method for value iteration. The method has been introduced in our previous work \cite{li2018approximate} and will be briefly reviewed here for completeness.

Given an \ac{mdp} $M = (S,A, P, s_0,\gamma, R)$, the objective is to find a policy $\pi$ that maximizes the total discounted return given by 
$J(s_0) = \max_\pi \Expect_\pi \sum_{t=0}^\infty \gamma^t R(s_t,a_t),$
where $s_t,a_t$ is the $t$-th state and action in the chain induced by policy $\pi$. 

The \ac{adp} for value iteration is to solve the following optimization problem:
\begin{align}
\label{eq:adp}
\min_\theta  & \sum_{s \in S} c(s;\theta)V(s;  \theta), \nonumber \\
\mbox{subject to: } & \calB V(s;\theta) - V(s; \theta) \leq 0 , \quad \forall s\in S,     
\end{align}
where the state relevant weight $c(s) = c(s; \theta)$ is the frequency with which different states are expected to be visited in the chain under policy $\pi(\cdot;\theta)$, and $ \calB V(s;\theta) = \tau \log \sum_{a} \exp \{ (R(s, a) + \gamma \sum_{s'}P(s' \mid s,a) V(s' ;\theta) / \tau \} $ is the softmax Bellman operator with the temperature $\tau >0$. The policy $\pi(\cdot;\theta):S\rightarrow \Delta(A)$ is computed from $V(\cdot;\theta)$ using \eqref{eq:policy_softmax}. It is shown in \cite{li2018approximate} that a value function $V(\cdot;\theta)$ satisfying the constraint in \eqref{eq:adp} is an upper bound on the optimal value function. The objective is to minimize this upper bound to approximate the optimal value function.

We  introduce a continuous function $B: \reals \rightarrow \reals_{+}$  with support equal to $[0, \infty)$. One such function is $B(x) = \max\{x,0\}$. Let $ g(s; \theta) =\calB V(s; \theta)- V(s; \theta), \forall s\in S$. Using randomized optimization~\cite{tadic2006randomized}, an equivalent representation of the set of constraints is	\begin{equation}
\label{eq:constraint_adp}
\Expect \limits_{ s \sim \Delta} B(g(s; \theta))=0,
\end{equation}
where $s$ is a random variable with a distribution $\Delta$ whose support is $S$. The augmented Lagrangian function of \eqref{eq:adp} with constraints replaced by \eqref{eq:constraint_adp} is
\begin{align}
\begin{split}
& L_{\nu}(\theta, \lambda, \xi) = \sum_{s\in S} c(s; \theta) V(s; \theta) + \lambda \cdot \Expect\limits_{s \sim \Delta} B(g(s; \theta))      \\
& + \frac{\nu}{2} \cdot \abs{\Expect\limits_{s \sim \Delta}B(g(s; \theta))}^2
\end{split}
\label{eq:rand-opt}
\end{align}
where $\lambda$ and $\xi$ are the Lagrange multipliers, and $ \nu$ is a large penalty constant. Using the Quadratic Penalty Function method  \cite{bertsekas1999nonlinear}, the solution is found by solving a sequence of optimization problems of the form: 
\begin{equation}	\label{eq:inner}
\min_{\theta \in \reals^\calK} \;    L_{\nu^k}(\theta, \lambda^k, \xi^k),
\end{equation}
where $\{\lambda^k\}$ and $\{\xi^k\}$ are sequences in $\reals$, $\{\nu^k\}$ is a positive penalty parameter sequence, and $\calK$ is the size of $\theta$.  After the inner optimization for \eqref{eq:inner} converges,  we update formula for multipliers $\lambda$ and $\xi$ in the \emph{outer optimization} as
\begin{align}
\label{eq:multipliers_update}
& \lambda^{k+1} = \lambda^{k} + \nu^{k} \cdot  \Expect\limits_{s \sim \Delta} B(g(s;\theta^k))
\end{align}
The outer optimization stops when it reaches the maximum number of iterations or $\norm{\nabla_{\theta}L_{\nu^k}(\theta^k, \lambda^k)} \leq \epsilon^{k}$. See \cite[Chap.~5.2]{bertsekas1999nonlinear} for more details about the augmented Lagrangian method.

By letting $c(s;\theta) =\sum_{t=0}^\infty P(X_t=s)$ be the state visitation frequency in the Markov chain $M^\theta$, for an arbitrary function $f:S\rightarrow \reals$,  it holds that $\underset{s \in S}{\sum}c(s; \theta)f(s;\theta) = \int p(h; \theta)f(h;\theta)dh$, where $p(h; \theta)$ is the probability of path $h$ in the Markov chain $M^\theta$, $f(h;\theta) = \sum_{i=1}^{|h|}f(s_{i};\theta)$.

By selecting $\Delta \propto c(\cdot;\theta)$ and letting $f(s;\theta) = V(s; \theta) + \lambda^k \cdot B( g(s; \theta)) + \frac{\nu^k}{2} \cdot \abs{B(g(s; \theta))}^2$, the $k$-th objective function  in \eqref{eq:inner} becomes $\min_{\theta} \underbrace{ \int p(h; \theta)f(h)dh}_ {F(\theta)}$. Parameter $\theta$ is updated by $\theta^{j+1} \leftarrow \theta^{j}  - \eta \cdot \nabla_{\theta} F(\theta)$, where $j$ represents the $j$-th inner iteration, $\eta$ is a positive step size. $\nabla_\theta F(\theta)  = \int \underbrace{\nabla_\theta p(h; \theta) f(h; \theta) dh}_{1} \\ + \int \underbrace{ p(h; \theta) \nabla_\theta f(h; \theta)dh}_{2}$, where using Monte-Carlo approximation, we have $1 \approx  \frac{1}{N_h} \sum_{h \sim p(h; \theta)} \bigl[\sum_{t=0}^{\abs{h}} \nabla_\theta \log \pi(a_t \mid s_t; \theta)\bigr] f(h; \theta), 2 \approx \frac{1}{N_h} \sum_{h \sim p(s_t; \theta)} \bigl[\sum_{t = 0}^{\abs{h}} \nabla_\theta f(s_t; \theta) \bigr]$, where $ \nabla_\theta f(s_t; \theta)= \nabla_{\theta} V(s_t; \theta) +\lambda^k \cdot \nabla_{g}B(g(s_t; \theta)) \nabla_{\theta}g(s_t; \theta) + \nu^k \cdot B(g(s_t; \theta) \nabla_{\theta}g(s_t; \theta)$.	

As the gradient of the objective function of the inner optimization problem can be computed from sampled trajectories, we can update the parameter $\theta$ using sampling-based augmented Lagrangian method, and thus have a model-free \ac{adp} method. The reader is referred to \cite{li2018approximate} for more technical details regarding the derivation of the method.

\bibliographystyle{ieeetr}
\bibliography{refs.bib}

\end{document}